\documentclass[reqno]{amsart}
\usepackage{graphicx}
\usepackage{amssymb}
\usepackage{amsfonts}
\usepackage{indentfirst}

\renewcommand\bigskip{\medskip}
\def\@oddhead{\hbox to \textwidth{\footnotesize{\it
Khintchine exponents and Lyapunov exponents } \hfill\thepage}}

\def\bc{\begin{center}}
\def\ec{\end{center}}

\newtheorem{thm}{Theorem}[section]

\newtheorem{lem}[thm]{Lemma}
\newtheorem{prop}[thm]{Proposition}

\numberwithin{equation}{section}

\newcommand{\set}[1]{\left\{#1\right\}}

\begin{document}
\title[On the frequency of partial quotients of regular continued fractions]
{On the frequency of partial quotients of regular continued fractions}

\author{Ai-Hua FAN}
\address{Ai-Hua FAN: \
Department of Mathematics,
 Wuhan University,
 Wuhan, 430072, P.R. China
 \& CNRS UMR 6140-LAMFA,
 Universit\'e de Picardie
 80039 Amiens, France}%
\email{ai-hua.fan@u-picaride.fr}%
\author{Ling-Min LIAO}
\address{Ling-Min LIAO: Department of Mathematics,
 Wuhan University,
 Wuhan, 430072, P.R. China \&
 CNRS UMR 6140-LAMFA,
 Universit\'e de Picardie
 80039 Amiens, France}%
 \email{lingmin.liao@u-picardie.fr}%
\author{Ji-Hua MA}
\address{Ji-Hua MA: Department of Mathematics,
 Wuhan University,
 Wuhan, 430072, P.R. China}%
 \email{jhma@whu.edu.cn}%





\maketitle

\begin{abstract}
We consider sets of real numbers in $[0,1)$ with prescribed
frequencies of partial quotients in their regular continued fraction
expansions. It is shown that the Hausdorff dimensions of these sets,
always bounded from below by $1/2$, are given by a modified
variational principle.
\end{abstract}


\section{Introduction}
Let $\mathbb{Q}^c$ denote the set of irrational number. It is
well-known that each $x\in [0,1)\cap\mathbb{Q}^c$ possesses a unique
 continued fraction expansion of the form
\begin{equation}\label{CF}
x=\frac{\displaystyle 1}{\displaystyle a_1(x)+ \frac{\displaystyle
1}{ a_2(x)+\displaystyle \frac{ 1}{ a_3(x)+\ddots}}},
\end{equation}
where $a_k(x)\in \mathbb{N}:=\set{1,2,3,\cdots}$ is the $k$-th {\em
partial quotient} of $x$. This expansion is usually denoted by $
x=[a_1(x),a_2(x), \cdots]$. For each $j\in\mathbb{N}$, define the
{\it frequency} of the digit $j$ in the continued fraction expansion
of $x$ by
\[\tau_j(x):= \lim_{n\to\infty} \frac{\tau_j(x,n)}{n},\]
whenever the limit exists, where $$\tau_j(x,n):= {\rm Card} \{1\le
k\le n: a_k(x)=j\}.$$

This paper is concerned with sets of real numbers with prescribed
digit frequencies in their continued fraction expansions. To be
precise, let $\vec{p}=(p_1,p_2,\dots)$ be a probability vector
with $p_j\geq 0$ for all $j\in \mathbb{N}$ and
$\sum_{j=1}^{\infty} p_j =1$, which will be called a {\it
frequency vector} in the sequel. Our purpose is to determine the
Hausdorff dimension of the set
\[ \mathcal{E}_{\vec{p}} := \{ x\in [0,1)\cap \mathbb Q ^c: \tau_j(x) =p_j \ \forall j\geq 1\}. \]

Let us first recall some notation. For any $a_1, a_2,\cdots, a_n\in
\mathbb{N}$, we call
\begin{equation*} I(a_1, a_2, \cdots, a_n):=\{x\in [0,1): a_1(x)=a_1,
a_2(x)=a_2, \cdots, a_n(x)=a_n\}
\end{equation*}
a {\em rank-$n$ basic interval}. Let $T:[0,1)\rightarrow [0,1)$ be
the Gauss transformation defined by
$$T(0)=0, \quad T(x)=1/x \ ({\rm mod}\ 1) \  \textrm{ for }x\in (0,1).$$
For a given frequency vector $\vec{p}=(p_1,p_2,\dots)$, we denote by
$\mathcal{N}(\vec{p})$ the set of $T$-invariant ergodic probability
measures $\mu$ such that
   \begin{equation}\label{measure}
\int |\log x| d \mu < \infty \textrm{ and }\mu(I(j))=p_j\textrm{
for all } j\geq 1.
   \end{equation}
Let $h_{\mu}$ stand for the measure-theoretical entropy of $\mu$,
and $\dim_H A$ for the Hausdorff dimension of a set $A$. The main
result of this paper can be stated as follows.
\begin{thm}\label{continuedfraction}
For any frequency vector $\vec{p} $, one has
\begin{align*}
  \dim_H(\mathcal{E}_{\vec{p}}) = \max \left\{ \frac{1}{2},
  \quad \sup_{\mu \in {\mathcal N}(\vec{p})} \frac{h_{\mu}}{2\int |\log x| d \mu } \right\},
  \end{align*}
\end{thm}
where the ``$\sup$" is set to be zero if ${\mathcal
N}(\vec{p})=\emptyset$.

  By virtue of $\log|T'(x)|=2|\log x|$, we see
that $2\int |\log x| d \mu$ is the Liapunov exponent of the measure
$\mu$. Therefore, the ``$\sup$" in the above is a variational
formula which relates the Hausdorff dimension to the entropy and
Liapunov exponent of measures.

  Theorem~\ref{continuedfraction} provides a complete solution to the
long standing problem that requests an exact formula for
$\dim_H(\mathcal{E}_{\vec{p}})$. Let us recall some partial results
in the literature. In 1966, under the condition that
$\sum_{j=1}^{\infty} p_j \log j < \infty$, Kinney and Pitcher
\cite{KP} showed that
\[ \dim_H (\mathcal{E}_{\vec{p}}) \geq \frac{-\sum_{j=1}^{\infty} p_j \log p_j}{2 \int |\log x| d
      \mu_{\vec{p}}},\]
where $\mu_{\vec{p}} $ is the Bernoulli measure on $[0,1]$ defined
by
$$\mu_{\vec{p}}(I(a_1, a_2, \cdots, a_n))=\prod_{j=1}^{n}p_{a_j}.$$
The above lower bound is just the Hausdorff dimension of the
Bernoulli measure $\mu_{\vec{p}} $. However, by a result of Kifer,
Peres and Weiss in 2001, this is not an optimal lower bound. Indeed,
it is shown in \cite{KPW} that, for any Bernoulli measure
$\mu_{\vec{p}}$,
\[ \dim_H \mu_{\vec{p}} \leq 1- 10^{-7}\]
This surprising fact indicates that the collection of Bernoulli
measures are insufficient for providing the correct lower bound for
$\dim_H(\mathcal{E}_{\vec{p}})$.

  In 1975, under the same condition
$\sum_{j=1}^{\infty} p_j \log j < \infty$, Billingsley and
Henningsen \cite{BH} obtained an improved lower bound
\[ \dim_H(\mathcal{E}_{\vec{p}})\geq  \sup_{\mu \in {\mathcal N}(\vec{p})} \frac{h_{\mu}}{2\int |\log x| d \mu }.\]
Moreover, they proved that, for any fixed $N\in\mathbb{N}$, this
lower bound is the exact Hausdorff dimension of the set
$$\set{x\in \mathcal{E}_{\vec{p}}:a_n(x)\leq N\textrm{ for all }n\geq 1 }$$
provided that $p_j=0$ for all $j>N$. It is therefore quite natural
to guess that this lower bound is the right value for $\dim_H
\mathcal{E}_{\vec{p}}$ in general. However, as will be shown in
Theorem \ref{continuedfraction}, this is not the case.

Actually, the lower bound due to Billingsley and Henningsen is only
a half of the correct lower bound. The other half of the lower
bound, namely, $\dim_H(\mathcal{E}_{\vec{p}})\geq 1/2$, can be
proved basing on Lemma 2.4 in \cite{LMW}. However, we will give a
direct proof in this paper.

The upper bound estimate is more difficult. In its proof, we will
use techniques from \cite{LMW} and \cite{BH} to  estimate  the
lengths of basic intervals. Not incidentally, an entropy-involved
combinatorial lemma borrowed from \cite{GW} (see Lemma
\ref{counting}) will play an important role.

The paper is organized as follows. In Section
\ref{Section-preliminary}, we give some preliminaries on the basic
intervals and on the entropy of finite words.
 In Section
\ref{Section-upper-bound}, we establish the upper bound in Theorem
\ref{continuedfraction}. In Section \ref{Section-lower-bound}, we
prove that $\dim_H(\mathcal{E}_{\vec{p}})\geq 1/2$ and show that we
can drop the condition $\sum_{j=1}^{\infty} p_j \log j < \infty$ in
Billingsley and Henningsen's theorem and then obtain the lower bound
in Theorem \ref{continuedfraction}. The last section serves as a
remark.

\medskip

\section{Preliminary}\label{Section-preliminary}

Let $x=[a_1(x), a_2(x), \cdots]\in [0,1)\cap\mathbb{Q}^c$.  The
$n$-th convergent in the continued fraction expansion of $x$ is
defined by
\begin{eqnarray*}
 \frac{p_n}{q_n}:=\frac{p_n(a_1(x),\cdots,a_n(x))}{q_n(a_1(x),\cdots,a_n(x))}
 =\frac{\displaystyle 1}{\displaystyle a_1(x)+
 \frac{\displaystyle 1}{\displaystyle a_2(x)+\ddots+\displaystyle
 \frac{1}{a_n(x)}}}.
 \end{eqnarray*}
For ease of notation, we shall drop the argument $x$ in what
follows. It is known (see \cite{Kh} p.9) that $p_n, q_n$ can be
obtained by the recursive relations:
\begin{eqnarray*}
p_{-1}=1, \ p_0=0, \ p_n = a_np_{n-1}+p_{n-2}  & & (n\geq 2),\\
q_{-1}=0, \ q_0=1, \ q_n = a_nq_{n-1}+q_{n-2}  & & (n\geq 2).
\end{eqnarray*}

By the above recursion relations, we have the following results.
\begin{lem}[\cite{Kh}]\label{q_n} Let $q_n=q_n(a_1,\cdots,a_n)$
and $p_n=p_n(a_1,\cdots,a_n)$, we have

\quad\quad \mbox{\rm (i)} \ $p_{n-1}q_n-p_nq_{n-1}=(-1)^n$;

\quad\quad \mbox{\rm (ii)} \ $q_n\geq 2^{\frac{n-1}{2}}$, \ \ \
$\prod\limits_{k=1}^na_k\leq q_n\leq \prod\limits_{k=1}^n(a_k+1).$
\end{lem}
\begin{lem}[\cite{Wu}]\label{delete-digits}For any $a_1, a_2,\cdots, a_n, b\in \mathbb{N}$,
\begin{eqnarray*} \frac{b+1}{2}\leq \frac{q_{n+1}(a_1,\cdots, a_j, b,
a_{j+1}, \cdots, a_n)}{q_n(a_1,\cdots, a_j, a_{j+1},\cdots,
a_n)}\leq b+1 \qquad (\forall 1\leq j <n).
\end{eqnarray*}\end{lem}

Recall that for any $a_1, a_2,\cdots, a_n\in \mathbb{N}$, the set
\begin{equation*} I(a_1, a_2, \cdots, a_n)=\{x\in [0,1): a_1(x)=a_1,
a_2(x)=a_2, \cdots, a_n(x)=a_n\}
\end{equation*}
is a rank-$n$ basic interval. We write $|I|$ for the length of an
interval $I$.

\begin{lem}[\cite{Kr1} p.18] \label{length}
The basic interval $I(a_1, a_2, \cdots, a_n)$ is an interval with
endpoints $p_n/q_n$ and $(p_n+p_{n-1})/(q_n+q_{n-1})$. Consequently,
one has
\begin{eqnarray}\label{length-basic-interval}
\Big|I(a_1,\cdots, a_n)\Big|=\frac{1}{q_n(q_n+q_{n-1})},
\end{eqnarray}
and
\begin{eqnarray}\label{length-estimate}
\frac{1}{2q^2_n}\leq \Big|I(a_1,\cdots, a_n)\Big| \leq
\frac{1}{q^2_n}.
\end{eqnarray}
\end{lem}

\begin{lem}\label{delete-estimate}
We have
\begin{eqnarray*}
|I(x_1,\cdots,j,\cdots, x_n)|\leq \frac{8}{(j+1)^2}
\Big|I(x_1,\cdots,\widehat{j},\cdots, x_n) \Big|,
\end{eqnarray*}
where the notation $\widehat{j}$ means ``deleting the digit $j$".
\end{lem}

\begin{proof}
   By
 Lemma \ref{delete-digits} and (\ref{length-estimate}), we have
\begin{eqnarray*}
|I(x_1,\cdots,j,\cdots, x_n)|&\leq&
{q^{-2}_n(x_1,\cdots,j,\cdots,x_n)} \\
&\leq& \left(\frac{2}{j+1}\right)^2
{q^{-2}_{n-1}(x_1,\cdots,\hat{j},\cdots, x_n)}\\
&\leq& \frac{8}{(j+1)^2} \Big|I(x_1,\cdots,\widehat{j},\cdots, x_n)
\Big|.
\end{eqnarray*}
\end{proof}

We will simply denote by $I_n(x)$ the rank $n$ basic interval
containing  $x$. Suppose that $a_n:=a_n(x)\geq 2$ and consider
     $$I_n^{\prime}(x)=I(a_1,\cdots,a_{n-1},a_n-1)\quad\textrm{and}\quad
     I_n^{\prime\prime}(x)=I(a_1,\cdots,a_{n-1},a_n+1)$$
which are two rank $n$ basic intervals adjacent to $I_n(x)$. By the
recursive equation of $q_n$ and (\ref{length-basic-interval}), one
has the following lemma.
   \begin{lem}\label{adjacent}
Suppose that $a_n:=a_n(x)\geq 2$. Then the lengths of the adjacent
intervals $I_n^{\prime}(x)$ and $I_n^{\prime\prime}(x)$  are bounded
by $\frac{1}{3}\bigl|I_n(x)\bigr|$ from below and by
$3\bigl|I_n(x)\bigr|$ from above.
   \end{lem}

For any $x\in [0,1)\setminus \mathbb{Q}$ and any word $i_1\cdots i_k
\in \mathbb{N}^k, (k\geq 1) $, denote by $\tau_{i_1\cdots i_k}(x,n)
$ the number of $j$, $1\leq j \leq n $, for which
\[a_j(x) \cdots a_{j+k-1}(x)=i_1\cdots i_k.\]

For $N\in \mathbb{N}$, define $\Sigma_N:=\{1,\dots, N\} $. We
shall use the following estimate in \cite{BH}.
\begin{lem}[\cite{BH}]\label{estimate-interval}
Let $N\geq 1$ and $n\geq 1$. For any $x=[x_1,x_2,\cdots]\in
[0,1]\cap \mathbb{Q}^c$ with $ x_j \in \Sigma_N$ for $1\leq j \leq n
$. Then for any $k\geq 1$, we have
\begin{eqnarray}\label{estimate-2.3}
\quad \log |I_n(x)| \leq 2\sum_{i_1\cdots i_k \in \Sigma_N^k}
\tau_{i_1\cdots i_k}(x,n) \log
  \frac{p_k(i_1,\cdots,i_k)}{q_k(i_1,\cdots,i_k)}+8+\frac{8n}{2^k}.
\end{eqnarray}
\end{lem}

Now we turn to the key combinatorial lemma which will be used in the
upper bound estimation. Let $\phi :[0,1]\rightarrow \mathbb{R} $
denote the function
\[
   \phi(0)=0, \quad \text{and} \quad \phi(t)=-t\log t \quad \text{for} \quad 0<t\leq 1.
\]
For every word $\omega \in \Sigma_N^n$ of length $n$ and every word
$ u\in \Sigma_N^k$ of length $k$, denote by $p(u|\omega) $ the
frequency of appearances of $u$ in $\omega$, i.e.,
\[
   p(u|\omega) = \frac{\tau_u(\omega)}{n-k+1},
\]
where $\tau_u(\omega)$ denote the number of $j$, $1\leq j \leq n-k+1
$, for which
\[\omega_j\cdots \omega_{j+k-1}=u.\]
Define
\[ H_k(\omega):=\sum_{u\in \Sigma_N^k} \phi(p(u|\omega)).\]
 We have the following counting lemma.
\begin{lem}[\cite{GW}]\label{counting}
 For any $h>0, \epsilon >0$,  $k\in \mathbb{N} $, and for any $n\in \mathbb{N}$ large
 enough, we have
 \[
    {\rm{Card}} \{\omega \in \Sigma_N^n : H_k(\omega) \leq k h\}
    \leq \exp (n(h+\epsilon)).
 \]
\end{lem}

\medskip

\section{Upper bound}\label{Section-upper-bound}

%
%
\subsection{Some Lemmas}

Let $(p(i_1,\cdots, i_k))_{i_1\cdots i_k \in \mathbb{N}^k} $ be a
probability vector indexed by $ \mathbb{N}^k$. As usual, we denote
by $q_k(a_1,\cdots, a_k) $ the denominator of the $k$-th convergent
of a real number with leading continued fraction digits $a_1, \dots,
a_k $.

\begin{lem}\label{less-than-1}
 For each $k\in \mathbb{N}$ and each probability vector
$(p(i_1,\cdots, i_k))_{i_1\cdots i_k \in \mathbb{N}^k}$,
\begin{align*}
\sum_{i_1\cdots i_k \in \mathbb{N}^k} -p(i_1,\cdots, i_k) \log
p(i_1,\cdots,i_k)
  \leq \sum_{i_1\cdots i_k \in \mathbb{N}^k}
-p(i_1,\cdots, i_k) \log
  |I(i_1,\cdots,i_k)|.
\end{align*}
\end{lem}

\begin{proof}
Applying Jesen's inequality to the concave function $\log $, we have
\[ \sum_{i_1\cdots i_k \in \mathbb{N}^k} p(i_1,\cdots, i_k) \log \frac{|I(i_1,\cdots,i_k)|}{p(i_1,\cdots, i_k)}
\leq \log  \sum_{i_1\cdots i_k \in \mathbb{N}^k}
|I(i_1,\cdots,i_k)|=0.
\]
\end{proof}

\smallskip

\begin{lem}\label{dominator-infinite}
Let $\vec{p}=(p_1,p_2,\dots)$ be a probability vector  and
 $\vec{q}=(q_1,q_2,\dots)$ a positive vector. Suppose
$-\sum_{j=1}^{\infty} p_j \log {q_j}= \infty $ and
$\sum_{j=1}^{\infty} q^{s}_j < \infty $ for some positive number
$s$. Then
\begin{eqnarray*}
\limsup_{n\to \infty} \frac{-\sum_{j=1}^{n} p_j \log p_j
}{-\sum_{j=1}^{n} p_j \log q_j } \leq s.
\end{eqnarray*}
\end{lem}

\begin{proof} This is a consequence of the following inequality (see
\cite{Wa}, p.217): for nonnegative numbers $s_j$ $(1\leq j\leq m)$
such that $\sum_{j=1}^{m} s_j =1$ and any real numbers $t_j \ (1\leq
j \leq m)$, we have
\begin{eqnarray}\label{entropy-inequality}
 \sum_{j=1}^{m} s_j( t_j- \log s_j) \leq  \log (\sum_{j=1}^{m}
 e^{t_j}).
\end{eqnarray}

 Fix $n\geq 1$.  Let
$s_j=p_j $ for $1\leq j\leq n $ and $s_{n+1}=
\sum_{j=n+1}^{\infty} p_j$. Let $t_j = s \log q_j$ for $1\leq j
\leq n$ and $t_{n+1}=0$. Applying the above inequality
(\ref{entropy-inequality}) with $m=n+1$, we get
\begin{eqnarray*}
s \sum_{j=1}^{n} p_j \log q_j-\sum_{j=1}^{n} p_j \log p_j
-(\sum_{j=n+1}^{\infty} p_j) \log (\sum_{j=n+1}^{\infty} p_j) \leq
\log (1+\sum_{j=1}^{n} q_j^{s}).
\end{eqnarray*}
Consequently,
\begin{eqnarray*}
 \frac{-\sum_{j=1}^{n} p_j \log p_j}{-\sum_{j=1}^{n} p_j \log
q_j}  \leq  s + \frac{(\sum_{j=n+1}^{\infty} p_j) \log
(\sum_{j=n+1}^{\infty} p_j)}{-\sum_{j=1}^{n} p_j \log q_j} +
\frac{\log (1+\sum_{j=1}^{n} q_j^{s})}{-\sum_{j=1}^{n} p_j \log
q_j}.
\end{eqnarray*}
Using the facts $-\sum_{j=1}^{\infty} p_j \log q_j = \infty$ and $
\sum_{j=1}^{\infty} q_j^{s}< \infty$, we  finish the proof by
letting $n\to\infty$.
\end{proof}

\smallskip
Lemma \ref{dominator-infinite} implies the following lemma. Recall
that $\Sigma_N^k= \{1,\dots, N\}^k$.
\begin{lem}\label{dominator-infinite-CF}
Let $k\geq 1$. If $(p(i_1,\cdots, i_k))_{i_1\cdots i_k \in
\mathbb{N}^k} $ is a probability vector such that \[\sum_{i_1\cdots i_k
\in \mathbb{N}^k} p(i_1,\cdots, i_k) \log q_k(i_1,\cdots,i_k)=
\infty,\] then we have
\begin{align*}
  \limsup_{N\to \infty} \frac{-\sum_{i_1\cdots i_k \in \Sigma_N^k} p(i_1,\cdots, i_k) \log p(i_1,\cdots,i_k)}
  {2\sum_{i_1\cdots i_k \in \Sigma_N^k} p(i_1,\cdots, i_k) \log
  q_k(i_1,\cdots,i_k)} \leq  \frac{1}{2}.
\end{align*}
\end{lem}

\begin{proof} By Lemma \ref{q_n}, for any $k\in \mathbb{N}$ and any $s> 1/2 $, we have
\begin{eqnarray*}
\sum_{i_1,\cdots, i_k} q_k(i_1,\cdots,i_k)^{-2s} \leq
\sum_{i_1,\cdots, i_k} (i_1\cdots
i_k)^{-2s}=(\sum_{j=1}^{\infty}j^{-2s})^k < \infty.
\end{eqnarray*}
 Thus
 we get the result by Lemma \ref{dominator-infinite}.
\end{proof}

\subsection{Proof of the upper bound}

\medskip

To prove the upper bound, we shall make use of  multi-step
Markov measures. Let $k\geq 1$, by a $(k-1)$-step Markov measure,
we mean a $T$-invariant probability measure $P$ on $[0,1)$
satisfying the Markov property
   \begin{equation}\label{Markov}
\frac{P(I(a_1,\cdots,a_n))}{P(I(a_1,\cdots,a_{n-1}))}
=\frac{P(I(a_{n-k},\cdots,a_n))}{P(I(a_{n-k},\cdots,a_{n-1}))}
   \end{equation}
for all $n\geq 1$ and $a_1,\cdots,a_n\in\mathbb{N}$ (see \cite{Ca},
p.9). We may regard a Bernoulli measure as a $0$-step Markov
measure.

For each $N\geq 2$, we denote by
$\mathbb{P}^k_{N}=\mathbb{P}^k_{N}(\vec{p})$ the collection of
$(k-1)$-step Markov measures satisfying the condition
   \begin{equation}\label{Markov-measure}
P(I(j))=p_j\textrm{ for }1\leq j\leq N-1\textrm{ and }
P(I(N))=1-\sum_{j=1}^{N-1}p_j.
   \end{equation}
These Markov measures are supported by the set of continued
fractions for which the partial quotients are bounded from above by
$N$.

 For each $i_1\cdots i_k \in
\set{1,2,\cdots,N}^k $, write $p(i_1,\dots,i_k)=P(I(i_1,\dots,
i_k))$. Put
\begin{eqnarray}\label{beta_N_k}
  \alpha_{N,k}:=  \sup_{P\in
  \mathbb{P}^k_{N}} \frac{-\frac{1}{k}\sum p(i_1,\cdots,i_k) \log p(i_1,\cdots,i_k)}{-2\sum p(i_1,\cdots,i_k) \log
  (p_k(i_1,\cdots,i_k)/q_k(i_1,\cdots,i_k))}.
\end{eqnarray}
The argument in \cite{BH} (pp.171-172) shows that the following
limit
\[
  {\alpha}_N:= \lim_{k\to \infty} \alpha_{N,k},
\]
exists and coincides with each of the following three limits:
\[
{\alpha}'_N:= \lim_{k\to \infty} \sup_{P\in
  \mathbb{P}^k_{N}} \frac{-\sum p(i_1,\cdots,i_k) \log p(i_1,\cdots,i_k)}{2\sum p(i_1,\cdots,i_k) \log
  q_k(i_1,\cdots,i_k)},
\]
\[
{\alpha}''_N:= \lim_{k\to \infty} \sup_{P\in
  \mathbb{P}^k_{N}} \frac{-\sum p(i_1,\cdots,i_k) \log p(i_1,\cdots,i_k)}{-\sum p(i_1,\cdots,i_k) \log
  |I(i_1,\cdots,i_k)|},
\]
and
\[
{\alpha}'''_N:= \lim_{k\to \infty} \sup_{P\in
  \mathbb{P}^k_{N}} \frac{h_P}{2\int |\log x| d P}.
\]

Let
\begin{eqnarray}\label{beta}
  \alpha:= \limsup_{N\to \infty} {\alpha}_N= \limsup_{N\to \infty}
{\alpha}'_N= \limsup_{N\to \infty} {\alpha}''_N=\limsup_{N\to
\infty} {\alpha}'''_N.
\end{eqnarray}
To prove the upper bound, we need only to prove the following two
propositions.

\begin{prop}\label{less-than-beta_N}
  For any $N\in \mathbb{N}$ large enough, we have
\begin{eqnarray*}
  \dim_H({\mathcal E}_{\vec{p}}) \leq \max \left\{\frac{1}{2}, \ \
 {\alpha}_N
  \right\}.
\end{eqnarray*}
\end{prop}
\begin{prop}\label{beta-less-than}
We have
  \begin{eqnarray*}
    {\alpha} \leq \max \left\{ \frac{1}{2},
  \quad \sup_{\mu \in {\mathcal N}(\vec{p})} \frac{h_{\mu}}{2\int |\log x| d \mu } \right\}.
\end{eqnarray*}
\end{prop}

Remark that we will finally establish the formula in Theorem
\ref{continuedfraction}, thus by Proposition \ref{less-than-beta_N}
and Proposition \ref{beta-less-than} we have
\[
\max \left\{\frac{1}{2}, \  \liminf_{N\to \infty}{\alpha}_N
\right\}=\max \left\{\frac{1}{2}, \  \limsup_{N\to \infty}{\alpha}_N
\right\}=\max \left\{ \frac{1}{2},
  \  \sup_{\mu \in {\mathcal N}(\vec{p})} \frac{h_{\mu}}{2\int |\log x| d \mu } \right\},\]
and if $\liminf_{N\to \infty}{\alpha}_N\geq 1/2 $, then the limit of
$\alpha_N $ exists and equals to
\[
\sup_{\mu \in {\mathcal N}(\vec{p})} \frac{h_{\mu}}{2\int |\log x| d
\mu }.
\]

\begin{proof}[Proof of Proposition~\textup{\ref{beta-less-than}.}]
By virtue of (\ref{beta}), we set
\[
\alpha = \limsup_{N\to \infty} \lim_{k\to \infty} \sup_{P\in
  \mathbb{P}^k_{N}} \frac{-\sum p(i_1,\cdots,i_k) \log p(i_1,\cdots,i_k)}{2\sum p(i_1,\cdots,i_k) \log q_k(i_1,\cdots,i_k)}.
\]
Denote
\[D:=\limsup_{N\to \infty} \lim_{k\to \infty} \sup_{P\in
  \mathbb{P}^k_{N}} {2\sum p(i_1,\cdots,i_k) \log q_k(i_1,\cdots,i_k)}.\]
If $D=\infty $, then by Lemma \ref{dominator-infinite-CF}, we have
$\alpha \leq 1/2 $.

Now suppose that $D<\infty $, which is equivalent to
\[\limsup_{N\to \infty} \lim_{k\to \infty} \sup_{P\in
  \mathbb{P}^k_{N}} {\int |\log x| d P}<\infty.\]
By (\ref{beta}),
\[
\alpha = \limsup_{N\to \infty} \lim_{k\to \infty} \sup_{P\in
  \mathbb{P}^k_{N}} \frac{h_P}{2\int |\log x| d P}.
\]
Without loss of generality, we may suppose that there is a sequence
of Markov measures $P_{N,k}\in \mathbb{P}^k_{N}$ converging to a
measure $\mu \in \mathcal{N}(\vec{p})$ in the $weak^{*}$-topology
such that
\[\limsup_{N\to \infty} \lim_{k\to \infty}
\frac{h_{P_{N,k}}}{2\int |\log x| d P_{N,k}}=\alpha.
\]
Then by the upper semi-continuity of the entropy function and the
weak convergence, we have
\[
{\alpha} \leq \sup_{\mu \in {\mathcal N}(\vec{p})}
\frac{h_{\mu}}{2\int |\log x| d \mu }.
\]
\end{proof}
%

\begin{proof}[Proof of Proposition~\textup{\ref{less-than-beta_N}.}]
For any fixed
integer $N$ which is large enough, and any $\epsilon>0$, we have
\begin{eqnarray*}
 E_{\vec{p}} \subset \bigcup_{\ell=1}^{\infty}
 \bigcap_{n=\ell}^{\infty}H_n(\epsilon,N),
\end{eqnarray*}
where
\[
H_n(\epsilon,N):=\left\{ x\in [0,1)\setminus \mathbb{Q}: \
\left|\frac{\tau_j(x,n)}{n}
 -p_j \right| < \epsilon, 1\leq j\leq N \right\}.
\]
For any $\gamma> \max \left\{1/2,
  {\alpha}_N
  \right\} $ and for any integer $k\in \mathbb{N}$, we have the
  $\gamma$-Hausdorff measure (see \cite{Fa1}, for the definition of ${\mathcal H}^{\gamma}$)
\begin{eqnarray*}
& &{\mathcal H}^{\gamma}\left(\bigcap_{n=\ell}^{\infty}H_n(\epsilon,N)\right)\\
 &\leq&
\sum_{|\frac{\tau_j(x,n)}{n}
 -p_j | < \epsilon, 1\leq j\leq N} |I_n(x)|^{\gamma} \quad \qquad (\forall n\geq \ell)\\
 &=& \sum_{n(p_j-\epsilon) < m_j <n(p_j+\epsilon), 1\leq j \leq N} \
 \sum_{x_1\cdots x_n \in A} |I(x_1,\cdots,x_n)|^{\gamma},
\end{eqnarray*}
where $$A:=\left\{x_1\cdots x_n \in \Sigma_N^n: \tau_j(x_1\cdots x_n)
= m_j, 1\leq j \leq N\right\}.$$
(We recall that
$\tau_j(x_1\cdots x_n)$ denotes the times of appearances  of $j$
in $x_1\cdots x_n$.)

 Let $\tilde{n}:=
\sum_{j=1}^N m_j $. By Lemma \ref{delete-estimate}, we have the
following estimate by deleting the digits $j>N$ in the first $n$
partial quotients $x_1,\dots,x_n $ of $x\in I_n(x)$:
\[
\sum_{x_1\cdots x_n \in A} |I(x_1,\cdots,x_n)|^{\gamma} \leq
\left(\sum_{j=N+1}^{\infty}\frac{8}{(j+1)^{2\gamma}}\right)^{n-\tilde{n}}
\sum_{x_1\cdots x_{\tilde{n}} \in \tilde{A}}
|I(x_1,\cdots,x_{\tilde{n}})|^{\gamma},
\]
where $$\tilde{A}:=\left\{x_1\cdots x_{\tilde{n}}\in
\Sigma_N^{\tilde{n}}: \tau_j(x_1\cdots x_{\tilde{n}}) = m_j, 1\leq j
\leq N\right\}.$$
 Since $\gamma>1/2$, for $N$ large enough
\begin{eqnarray}\label{part1}
    \sum_{j=N+1}^{\infty}\frac{8}{(j+1)^{2\gamma}} <1.
\end{eqnarray}

 By applying Lemma \ref{estimate-interval}, and
noticing that $\tau_{i_1\cdots i_k}(x,{\tilde{n}}) \leq
\tau_{i_1\cdots i_k}(x_1\cdots x_{\tilde{n}})+k $, we have
\begin{eqnarray*}
& &|I(x_1,\cdots,x_{\tilde{n}})|= \exp \{\log
|I(x_1,\cdots,x_{\tilde{n}})| \}\\
 &\leq&  \exp \left\{2\sum_{i_1\cdots i_k \in
\Sigma_N^k} (\tau_{i_1\cdots i_k}(x_1\cdots x_{\tilde{n}})+k) \log
  \frac{p_k(i_1,\cdots,i_k)}{q_k(i_1,\cdots,i_k)}+8+\frac{8{\tilde{n}}}{2^k}\right\}.
\end{eqnarray*}
Thus
\begin{eqnarray*}
 & &\sum_{x_1\cdots x_{\tilde{n}} \in \tilde{A}}
|I(x_1,\cdots,x_{\tilde{n}})|^{\gamma}\\
 &\leq& \sum_{m_{i_1 \cdots i_k}} \sum_{x_1\cdots x_{\tilde{n}}\in B}
\exp \left\{2\gamma \sum_{i_1\cdots i_k \in \Sigma_N^{k}}
(m_{i_1\cdots i_k}+k) \log
  \frac{p_k}{q_k}+8\gamma+\frac{8{\tilde{n}}\gamma}{2^k}\right\},
\end{eqnarray*}
where $$B:=\left\{x_1\cdots x_{\tilde{n}}\in \tilde{A}:
{\tau_{i_1\cdots i_k}(x_1\cdots x_{\tilde{n}})
 =m_{i_1\cdots i_k}} \ \forall i_1\cdots i_k \in
 \Sigma_N^{k}\right\}.
 $$

Take \begin{eqnarray}\label{h} h= \frac{1}{k} \sum_{i_1\cdots i_k
\in \Sigma_{N}^k} \phi\left(\frac{m_{i_1\cdots
i_k}}{\tilde{n}-k+1}\right)
\end{eqnarray}
in Lemma \ref{counting}. We have for any $\delta>0$ and for
$\tilde{n}$ large enough
\begin{eqnarray*}
& &\sum_{x_1\cdots x_{\tilde{n}}\in B} \exp \left\{2\gamma
\sum_{i_1\cdots i_k \in \Sigma_N^{k}} (m_{i_1\cdots i_k}+k) \log
  \frac{p_k}{q_k}+8\gamma+\frac{8{\tilde{n}}\gamma}{2^k}\right\}\\
&\leq& \exp \left\{\tilde{n}\left(h+\delta \right)+ 2\gamma
\sum_{i_1\cdots i_k \in \Sigma_N^{k}} (m_{i_1\cdots i_k}+k) \log
  \frac{p_k}{q_k}+8\gamma+\frac{8{\tilde{n}}\gamma}{2^k}\right\}.
\end{eqnarray*}

Rewrite the right side of the above inequality as \[\exp
\left\{\tilde{n}\left(L(\gamma, k,m_{i_1\cdots
i_k})\right)\right\},\] where
\begin{eqnarray*}
L(\gamma, k,m_{i_1\cdots i_k}):=h+ 2\gamma \sum_{i_1\cdots i_k \in
\Sigma_N^{k}} \frac{m_{i_1\cdots i_k}+k}{\tilde{n}} \log
  \frac{p_k}{q_k}+\frac{8\gamma}{\tilde{n}}+\frac{8\gamma}{2^k} + \delta .
\end{eqnarray*}
Since there are at most $(\tilde{n}-k+1)^{N^k} $ possible words of
$i_1\cdots i_k $ in $\Sigma_N^{\tilde{n}}$, we have
\begin{eqnarray*}
 & &\sum_{x_1\cdots x_{\tilde{n}} \in \tilde{A}}
|I(x_1,\cdots,x_{\tilde{n}})|^{\gamma}\\
 &\leq& (\tilde{n}-k+1)^{N^k}
\exp \left\{\tilde{n} \left( \sup_{m_{i_1\cdots i_k}}L(\gamma,
k,m_{i_1\cdots i_k}) \right)\right\}.
\end{eqnarray*}

Notice that by the definition of $\tilde{A} $ and $B$, the possible
values of $m_{i_1\cdots i_k} $ are restricted to satisfy the
condition that the frequency of digit $j$ in $x_1\cdots
x_{\tilde{n}}$ is about $p_j$. Then when $\tilde{n} \to \infty$, for
${i_1,\cdots,i_k}\in \Sigma_N^k $
\begin{eqnarray}\label{appro-Markov}
\frac{ m_{i_1\cdots i_k}}{\tilde{n}-k+1} \rightarrow
p(i_1,\cdots,i_k),
\end{eqnarray}
and $\{p(i_1,\cdots,i_k): {i_1,\cdots,i_k}\in \Sigma_N^k \} $
defines a probability measure $P$ in $\mathbb{P}_N^k $.

Now take $\delta>0 $ small enough and $k$ large enough such that
\begin{eqnarray}\label{delta1}
  \frac{8\gamma}{2^k}< \delta,
\end{eqnarray}
and
\begin{eqnarray}\label{delta-gamma}
\gamma> \frac{-\frac{1}{k}\sum p(i_1,\cdots,i_k) \log
p(i_1,\cdots,i_k)+ 5 \delta}{-2\sum p(i_1,\cdots,i_k) \log
  (p_k(i_1,\cdots,i_k)/q_k(i_1,\cdots,i_k))}.
\end{eqnarray}
The last inequality comes from the definition of $\alpha_N$,
(\ref{beta_N_k}) and the assumption $\gamma > \alpha_N $.

By (\ref{appro-Markov}), for sufficiently large $\tilde{n}$, we have
$\frac{8\gamma}{\tilde{n}}< \delta$ and
\begin{eqnarray}\label{delta2}
\begin{split}
\left|\frac{1}{k} \sum_{i_1\cdots i_k \in \Sigma_N^k} \phi\left(
\frac{m_{i_1\cdots i_k}}{\tilde{n}-k-1}\right) -\frac{1}{k}
\sum_{i_1\cdots i_k \in \Sigma_N^k} \phi\left( p(i_1\cdots
i_k)\right) \right|< \delta, \\
\left| \sum_{i_1\cdots i_k \in \Sigma_N^k} \frac{m_{i_1\cdots
i_k}+k}{\tilde{n}}\log \frac{p_k}{q_k} - \sum_{i_1\cdots i_k \in
\Sigma_N^k} p(i_1\cdots i_k)\log \frac{p_k}{q_k} \right|< \delta.
\end{split}
\end{eqnarray}
By (\ref{h}) and (\ref{delta2}),
\[L(\gamma,k,m_{i_1\cdots i_k})<\frac{1}{k}
\sum_{i_1\cdots i_k \in \Sigma_N^k} \phi\left( p(i_1\cdots
i_k)\right)+2\gamma \sum_{i_1\cdots i_k \in \Sigma_N^k} p(i_1\cdots
i_k)\log \frac{p_k}{q_k}+5\delta. \] Thus (\ref{delta-gamma})
implies
\[
  L(\gamma,k,m_{i_1\cdots i_k}) <0.
\]
 Hence finally, by
(\ref{part1}), we can obtain for any $\gamma> \max \left\{1/2,
\alpha_N \right\} $,
\[{\mathcal H}^{\gamma}\left(\bigcap_{n=k}^{\infty}H_n(\epsilon,N)\right) < \infty.\]

This implies that $\dim_H({\mathcal E}_{\vec{p}})\leq \max
\left\{1/2, \alpha_N \right\}$ as desired.
\end{proof}

\medskip

\section{Lower bound}\label{Section-lower-bound}
In this section, we first prove $\dim_H (\mathcal E_{\vec{p}}) \geq
{1}/{2}$. Then we examine what happens if the  condition
$\sum_{j=1}^{\infty} p_j \log j < \infty$ in Billingsley and
Henningsen's theorem is violated. We will see that if the condition
is not satisfied, then $\dim_H (\mathcal E_{\vec{p}}) = {1}/{2}$ and
$\mathcal{N}(p)=\emptyset$.

The following is the key lemma for proving that $\dim_H (\mathcal
E_{\vec{p}}) \geq {1}/{2}$.
\begin{lem}\label{lem-seed}
For any given sequence of positive integers $\{c_n \}_{n\geq 1} $
tending to the infinity, there exists a point $z=(z_1,z_2, \dots)
\in \mathcal{E}_{\vec{p}}$ such that $ z_n \leq c_n$ for all $n\geq
1$.
\end{lem}
\begin{proof} For any $n\geq 1$, we construct a probability vector
$(p_1^{(n)}, p_2^{(n)}, \dots, p_{k}^{(n)},\dots) $ such that
$p_k^{(n)}
>0 $ for all $1\leq k \leq c_n $ and $\sum_{k=1}^{c_n} p_k^{(n)} =1 $, and that for any $k\geq
1 $,
\begin{eqnarray}\label{p-lim}
\lim_{n\to \infty} p_k^{(n)}=p_k.
\end{eqnarray}
Consider a product Bernoulli probability $\mathbb{P}$ supported by
$\prod_{n=1}^{\infty} \{ 1, \dots, c_n \}$. For each digit $k\geq
1$, consider the random variables of $x\in
{\mathbb{N}}^{\mathbb{N}}$, $X_i(x)= {1}_{\{k\}} (x_i)$, $ ( i\geq
1)$. By Kolmogorov's strong law of large numbers (see \cite{Sh}
p.388), we have for each digit $k$,
\begin{eqnarray*}
 \lim_{n\to \infty} \frac{1}{n} \left(\sum_{i=1}^{n} {1}_{\{k\}} ( x_i)  -
 \sum_{i=1}^{n} \mathbb{E} ({1}_{\{k\}} ( x_i))\right) =0
 \qquad \mathbb{P}-a.s.,
\end{eqnarray*}
which implies
\begin{eqnarray}\label{pk-frequency}
 \lim_{n\to \infty} \frac{1}{n} \sum_{i=1}^{n} {1}_{\{k\}}
 (x_i)
 = \lim_{n\to \infty} \frac{1}{n} \sum_{i=1}^{n} p_k^{(i)} =p_k
 \qquad \mathbb{P}-a.s..
\end{eqnarray}
 That is to say, for $\mathbb{P}$ almost every point in
the space $\prod_{n=1}^{\infty} \{ 1, \dots, c_n \} $, the digit $k$
has the frequency $p_k$. Considering each point in
$\mathbb{N}^{\mathbb{N}}$ as a continued fraction expansion of a
number in $[0,1] $, we complete the proof.\end{proof}

\begin{proof}[Proof of ~\textup{$\dim_H (\mathcal E_{\vec{p}}) \geq
{1}/{2}$.}]
Take $c_n=n$ in Lemma \ref{lem-seed}, we find a
point $z\in \mathcal E_{\vec{p}}$, such that
\begin{eqnarray}\label{z-n-less}
z_n = a_n(z) \leq n \qquad (\forall n \geq 1).
\end{eqnarray}

For a positive number $b >1$, set
 \[ \mathcal F _z(b):= \{ x\in [0,1):  a_{k^2}(x) \in
  (b^{k^2},2b^{k^2}]; \ \ a_{k}(x)=a_{k}(z) \ {\rm if} \ k \ {\rm is \ nonsquare}
    \}.
 \]
It is clear that $\mathcal F _z(b) \subset \mathcal E_{\vec{p}}$ for
all $b>1$. We define a measure $\mu$ on $\mathcal F_z(b)$. For $n^2
\leq m < (n+1)^2$, set
\begin{eqnarray}\label{measure-Luroth}
  \mu (I_m(x)) = \prod_{k=1}^n \frac{1}{b^{k^2}}.
\end{eqnarray}
Denote by $B(x,r)$ the ball centered at $x$ with radius $r$. We will
show that for any $\theta>0$, there exists $b>1$, such that for all
$x\in \mathcal F _z(b)$,
\begin{eqnarray}\label{local_density}
\liminf_{r\to 0} \frac{\log \mu(B(x,r))}{\log r} \geq
\frac{1}{2}-\theta.
\end{eqnarray}

In fact, for any positive number $r$, there exist integers $m$ and
$n$ such that
\begin{eqnarray}\label{length-r}
|I_{m+1}(x)| < 3r \leq |I_m(x)| \quad \text{and} \quad n^2\leq m <
(n+1)^2.
\end{eqnarray}

 By the construction
of $\mathcal F_z(b) $, $a_{n^2}(x)> b^{n^2} >1$. Let
$x=[x_1,x_2,\cdots] $. By Lemma \ref{adjacent}, $B(x,r)$ is covered
by the union of three adjacent rank $n^2$ basic intervals, i.e.,
\begin{eqnarray*}
  B(x,r) \subset I(x_1,x_2,\cdots,
x_{n^2}-1) \cup I(x_1,x_2,\cdots, x_{n^2}) \cup I(x_1,x_2,\cdots,
x_{n^2}+1).
\end{eqnarray*}
By the definition of $\mu$, the above three intervals admit the same
measure.
Hence by (\ref{length-r}), we have
\begin{eqnarray}\label{mu_B}
\frac{\log \mu(B(x,r))}{\log r} \geq \frac{ \log
3\mu(I(x_1,x_2,\cdots, x_{n^2}))}
 {\log \frac{1}{3}|I(x_1,x_2,\cdots, x_{m+1})|}.
\end{eqnarray}

However, on the one hand,  by (\ref{measure-Luroth})
\begin{eqnarray}\label{estimate-numerator}
 -\log \mu(I(x_1,x_2,\cdots, x_{n^2}))
=-\log \prod_{k=1}^{n} \frac{1}{b^{k^2}}= \sum_{k=1}^{n} k^2 \log b.
\end{eqnarray}

On the other hand, by (\ref{length-estimate}) and Lemma \ref{q_n},
we have
\begin{eqnarray*}
 -\log |I(x_1,x_2,\cdots, x_{m+1})| \leq \log 2+ \sum_{k=1}^{m+1} 2\log (x_k+1).
\end{eqnarray*}
Let us estimate the second term of the sum. First we have
\begin{eqnarray*}
  \sum_{k=1}^{m+1} 2\log (x_k+1)\leq  2\sum_{k=1}^{n+1} \log (x_{k^2}+1) + 2\sum_{k=1}^{m+1}
\log (z_k+1).
\end{eqnarray*}
Since $x_{k^2} \leq 2b^{k^2} $ for all $k\geq 1$, we deduce
\begin{eqnarray*}
  \sum_{k=1}^{m+1} \log (x_{k^2}+1)&\leq& \sum_{k=1}^{n+1} \log (2b^{k^2}+1)
   \leq  \sum_{k=1}^{n+1} \log (3b^{k^2})  \\
&= &  (n+1)\log 3+ \sum_{k=1}^{n+1} k^2 \log b.
\end{eqnarray*}
By (\ref{z-n-less}), since $z_n\leq n $, for all $n\geq 1 $, we know
\begin{eqnarray*}
  \sum_{k=1}^{m+1}
\log (z_k+1)\leq \sum_{k=1}^{(n+1)^2} \log (k+1).
\end{eqnarray*}
Thus
\begin{eqnarray}\label{estimate-denominator}
\begin{split}
&-\log |I(x_1,x_2,\cdots, x_{m+1})| \\
\leq& \log2 +
  2(n+1)\log3+2\sum_{k=1}^{n+1}k^2 \log b + 2\sum_{k=1}^{(n+1)^2}
  \log (k+1).
  \end{split}
\end{eqnarray}

Combining (\ref{estimate-numerator}) and
(\ref{estimate-denominator}), for any $\theta>0$, take $b>1$ to be
large enough, we have
\begin{eqnarray*}
\liminf_{n\to\infty} \frac{ \log \mu(I(x_1,x_2,\cdots, x_{n^2}))}
 {\log |I(x_1,x_2,\cdots, x_{m+1})|} \geq \frac{1}{2}- \theta,\quad \forall x\in \mathcal F_z(b).
\end{eqnarray*}
Hence by (\ref{mu_B}), we obtain (\ref{local_density}).


Since $\theta$ can be arbitrary small, by Billingsley Theorem
(\cite{Bi}), we have
\[
  \dim_H (\mathcal E_{\vec{p}}) \geq \frac{1}{2}.
\]
\end{proof}

\medskip
Thus Theorem \ref{continuedfraction} is already proved under  the
condition $\sum_{j=1}^{\infty} p_j \log j < \infty$.

Now assume $\sum_{j=1}^{\infty} p_j \log j = \infty$. Then for any
invariant and ergodic measure $\mu$ such that $ \mu(I(j))=p_j$ for
all $j\geq 1$, we have
\begin{eqnarray}\label{Lya-infinite}
 \int |\log x| d\mu(x) \geq \sum_{j=1}^{\infty} \mu(I(j))\log j= \sum_{j=1}^{\infty} p_j\log
 j=\infty.
\end{eqnarray}
This implies $D=\infty$ (see the proof of Proposition
\ref{beta-less-than} for the definition of $D$). Then by Proposition
\ref{less-than-beta_N}, we have $\dim_H (\mathcal E_{\vec{p}}) \leq
\frac{1}{2}$. Since we have already proved $\dim_H (\mathcal
E_{\vec{p}}) \geq \frac{1}{2} $, we get
\[
  \dim_H (\mathcal E_{\vec{p}}) = \frac{1}{2} \quad {\rm if} \quad \sum_{j=1}^{\infty} p_j \log j = \infty.
\]
This is in accordance with the formula of Theorem
\ref{continuedfraction} under the convention that $\sup \emptyset =0
$, because (\ref{Lya-infinite}) implies
$\mathcal{N}(\vec{p})=\emptyset$.

Finally, we remark that $\mathcal{N}(\vec{p})=\emptyset $ if and
only if $\sum_{j=1}^{\infty} p_j \log j = \infty$. We have seen the
``if" part. For the other part, assume $\sum_{j=1}^{\infty} p_j \log
j < \infty$. Then the Bernoulli measure $\mu$ such that
$\mu(I(j))=p_j$ satisfies
\begin{eqnarray*}
 \int |\log x| d\mu(x) \leq \sum_{j=1}^{\infty} \mu(I(j))\log (j+1)= \sum_{j=1}^{\infty} p_j\log
 (j+1)<\infty.
\end{eqnarray*}
which implies that $\mathcal{N}(\vec{p}) \neq \emptyset $.

\section{A remark}

As suggested by the referee, we add a remark on a problematic
argument appearing in the literature. To obtain an upper bound of
the Hausdorff dimension of a set, one usually applies the
Billingsley's theorem by constructing a finite measure $P$ on the
set such that
 \begin{eqnarray*}
       |U|^s \leq P(U)
     \end{eqnarray*}
(see \cite{Fa1}, p.67). For the set ${\mathcal{E}}_{\vec{p}}$, where
$p_j=0$ for some $j$, the Markov measure $P$ satisfying
(\ref{Markov-measure}) does not match because the cylinders starting
with $j$ do not charge the measure and the above inequality is
obviously not true. This appeared unnoticed for long (see the
remarks of Kifer \cite{Kifer}, p. 2012).

This problem did exist in the proof of Theorem 2 in \cite{BH}. Let
us briefly indicate how to get around the problem in the proof of
Theorem 2 of \cite{BH} when
  $p_j=0$ for some $j$'s. The basic idea is similar to that of
 Cajar (see \cite{Ca}, p.67) and that of Kifer \cite{Kifer}.
Recall that $P$ is a $(k-1)$-step Markov measure supported on the
set
  $$\set{x\in [0,1)\cap\mathbb{Q}^c:a_n(x)\leq N\textrm{ for all
}n\geq 1 }.$$ It is uniquely determined by its values on the
$k$-cylinders, namely,
$$p(i_1,\dots,i_k)=P([i_1,\dots, i_k]),\quad (i_1\cdots i_k)
\in \set{1,2,\cdots,N}^k.$$ Let $0<\epsilon<1$ and $P_{\epsilon}$ be
the perturbed $(k-1)$-Markov measure determined by
$$p_{\epsilon}(i_1,\dots,i_k)=(1-\epsilon)P([i_1,\dots, i_k])+\frac{\epsilon}{N^k},\quad (i_1\cdots i_k)
\in \set{1,2,\cdots,N}^k.$$ Now, we can apply the Billingsley's
theorem with $P_{\epsilon}$ to find an upper bound, and then get the
desired result
  by letting $\epsilon\to 0$.

In the present paper, we have intentionally avoided using the
Billingsley's theorem. The proof of upper bound consists of
Propositions \ref{less-than-beta_N} and \ref{beta-less-than}.
Proposition \ref{beta-less-than} concerns some calculations for
which the zero frequency of some digits will not cause any trouble.
In the proof of Proposition \ref{less-than-beta_N}, we have used a
``covering argument" depending on the estimate (\ref{estimate-2.3})
instead of using the Billingsley's theorem.  This enables us to get
the upper bound of the Hausdorff dimension.

\medskip

{\it{Acknowledgments : }} The authors are grateful to the referee
for a number of valuable comments. This work was partially supported
by NSFC10728104 (A. H. Fan) and NSFC10771164 (J. H. Ma).

\bibliographystyle{amsplain}
{ }

\end{document}